\documentclass[11pt]{amsart}
\usepackage{latexsym,amssymb,amsmath}
\def\ctimes{\times \cdots\times}

\def\s{\sigma}

\def\trank{\text{rank}}

\def\BC{\mathbb C}

\def\BP{\mathbb P}
\def\pp#1{\mathbb P^{#1}}

\def\pp#1{{\mathbb P}^{#1}}
\def\tdim{\rm dim}
\def\hd{,...,}

\def\cB{{\mathcal B}}

\def\cG{{\mathcal G}}\def\cM{{\mathcal M}}
\def\cR{{\mathcal R}}

\def\cC{{\mathcal C}}

\def\cO{{\mathcal O}}\def\cV{{\mathcal V}}
\def\cQ{\mathcal Q}

\def\11{\mathbf 1}

\def\FS{{\mathfrak S}}

\def\a{\alpha}

\def\s{\sigma}

\def\ot{{\mathord{\,\otimes }\,}}
\def\op{{\mathord{\,\oplus }\,}}
\def\otc{{\mathord{\otimes\cdots\otimes}\;}}

\def\ctimes{{\mathord{\times\cdots\times}\;}}

\def\ra{{\mathord{\;\rightarrow\;}}}

\def\La#1{\Lambda^{#1}}
\def\tim{\text{Image}\,}

\def\tdim{\text{dim}\,}
\def\tcodim{\text{codim}\,}

\def\trank{\text{rank}\,}

\def\tdet{\text{det}}

\newtheorem{theo}{Theorem}

\newtheorem{theorem}{Theorem}[section]
\newtheorem{proposition}[theorem]{Proposition}
\newtheorem{lemma}[theorem]{Lemma}

\theoremstyle{definition}
\newtheorem{definition}[theo]{Definition}

\newtheorem{example}[theo]{Example}

\theoremstyle{remark}
\newtheorem{remark}[theorem]{Remark}

\begin{document}

\title{On the ideals and singularities of secant varieties of Segre varieties}
\author{J.M. Landsberg 
 \and Jerzy Weyman}
\begin{abstract}We   find 
generators for the ideals of secant varieties
of Segre varieties in the cases
of $\s_k(\pp 1\times \pp n\times \pp m)$ for all $k,n,m$, 
$\s_2(\pp n\times\pp m\times \pp p\times \pp r)$ for all $n,m,p,r$ (GSS conjecture for four factors),
and $\s_3(\pp n\times\pp m\times \pp p)$ for all $n,m,p$ 
and prove they are normal with rational singularities
in the first case and
arithmetically Cohen-Macaulay in the second two.  
\end{abstract}
\thanks{Supported respectively by NSF grants DMS-0305829 and DMS-0600229}
\thanks{MSC 13P99,14Q15,15A69}
\email{jml@math.tamu.edu, j.weyman@neu.edu}

\maketitle

\section{Introduction}

Let $V$ be a vector space over a field $K$ of characteristic zero  and
let $X \subset\BP V$ be a projective variety. Define $\s_r(X)$,
 the {\it variety of
secant $\pp{r-1}$'s to $X$} by
$$
\s_r(X)=\overline{ \cup_{x_1\hd x_r\in X}\BP_{x_1\hd x_r }}
$$
where $\BP_{x_1\hd x_r }\subset\BP V$ denotes the linear space
spanned by $x_1\hd x_r$ (usually a $\BP^{r-1}$) and the overline denotes
  Zariski closure.

Let $A_1\hd A_n$ be vector spaces over   $K$,
with $\tdim A_j=a_j$.
Let $Seg(\BP A_1^*\ctimes \BP A_n^*)\subset \BP (A_1^*\otc A_n^*)$
denote the {\it Segre variety} of decomposable tensors.
(We use the dual vector spaces $A_j^*$ when discussing varieties
because we will mostly be concerned with modules of polynomials
and this convention enables our modules to be $*$-free.)

For  applications
to computational complexity, algebraic statistics,
and other areas, one would like to have the defining
equations for secant varieties of Segre varieties
$\s_r(\BP A_1^*\ctimes \BP A_n^*)=\s_r(Seg(\BP A_1^*\ctimes \BP A_n^*))$
and understand their singularities. In computational
complexity  one studies the stratification
of $A_1^*\ot A_2^*\ot A_3^*$ by the secant varieties
of the Segre,
as given a bilinear map $f: A_1 \times A_2 \ra A_3^*$
(such as matrix multiplication when each
$A_j$ is the space of $m\times m$ matrices), the
smallest $r$ such that $f\in \s_r( \BP A_1^*\times \BP A_2^*\times \BP A_3^*)$
is a measure of its complexity.
More generally, in algebraic statistics (see, e.g,
\cite{berndetal}), one
would like as much information
as possible about different   algebraic
statistical models, and secant varieties of
Segre varieties are important special classes of such models.
The techniques employed in this paper will be
useful for the general study of these models.

\smallskip 
 
Remarkably little is known about even
set-theoretic defining equations of the
$\s_r(\BP A_1^*\ctimes \BP A_n^*)$, let alone generators
of the ideals (which is considerably more difficult).
The only case well understood is the case $n=2$ where the secant varieties are the classical determinantal varieties. In the case   $n=3$, the defining ideal of   $\s_2(\BP A_1^*\times \BP A_2^*\times \BP A_3^*)$ was described in  \cite{LMsec}, using the
methods of \cite{weyman}.  
Set theoretic generators for $\s_2(\BP A_1^*\times \cdots\times \BP A_n^*)$ were also described in \cite{LMsec}.

In the present paper we take the next step in
 understanding generators of the   ideals and 
 singularities  of the varieties
$\s_r(\BP A_1^*\times \cdots\times \BP A_n^*)$.
We make extensive use of the machinery of \cite{weyman}.

A significant role in our study is played by 
   auxiliary varieties that contain
$\s_r(\BP A_1^*\times \cdots\times \BP A_n^*)$  and
have ideals that are easier to study. The simplest
of these is the following:

\begin{definition}
Let $b_j\leq a_j:=\tdim A_j$ be nonnegative integers. Define
the {\it subspace varieties}
$$
Sub_{b_1\hd b_n}:=
\{ T\in A_1^*\otc A_n^* \mid \exists A_j^{'*}\subseteq A_j^*,\
\tdim A_j'=b_j,\ T\in A_1^{'*}\otc A_n^{'*}\}.
$$
\end{definition}

Subspace varieties are cousins of the {\it rank varieties} in \cite{weyman}.
We use the terminology \lq\lq subspace\rq\rq\  to avoid confusion with tensor
rank. They are useful because $ \s_r(\BP A_1^*\ctimes \BP A_n^*)
\subset \BP(Sub_{r\hd r})$ and $Sub_{r\hd r}$ admits
a nice desingularization described in \S\ref{subpfsect}.

We first determine generators of the ideals of
the subspace varieties using elementary representation
theory and prove that they
are normal, with rational singularities using techniques
from \cite{weyman}
in \S\ref{subpfsect}.  Subspace varieties enable one to
reduce the problem of finding generators
of the ideals of
$\s_r(\BP A_1\times \cdots\times \BP A_n)$ where
$\tdim A_j\geq r$ to the cases where
$\tdim A_j=r$ for all $j$ (Proposition \ref{gensofsr}),
which we refer to as the {\it basic cases}.
In \S\ref{3p2sect} and \S\ref{4p1sect} we respectively resolve
the basic cases of $\s_3(\pp 2\times\pp 2\times \pp 2)$
and $\s_2(\pp 1\times\pp 1\times \pp 1\times \pp 1)$.

Recall that for any variety $Z\subset \BP V$ invariant
under the action of an algebraic group $G$, the
generators of the ideal of $Z$ will be
grouped into $G$-modules. In our
case $G=SL(A_1)\times\cdots\times SL(A_n)$, and
the special linear group has the added feature that
the decomposition of its various modules is
essentially independent of the dimension of the
vector space $A_j$. For example, when $n=2$,  the ideal of
$ \s_r(\BP A^*\times \BP B^*) $ is generated
by the irreducible module $\La{r+1}A\ot \La{r+1}B
\subset S^{r+1}(A\ot B)$ which corresponds to the
space of
$r+1\times r+1$ minors as long as
$\tdim A, \tdim B\geq r+1$.

Finally, recall that a {\it flattening}
of a tensor $T\in A_1\otc A_n$ is
to let $I=\{i_1\hd i_p\}\subset \{1\hd n\}$, $J=\{1\hd n\}\backslash I$,
$A_I=A_{i_1}\otc A_{i_p}$, $A_J=A_{j_1}\otc A_{j_{n-p}}$
and consider $T\in A_I\ot A_J$.
Flattenings are useful because
  the
ideals of secant varieties of Segre products
of two projective spaces are well understood.

\subsection*{Notation} For a partition $\pi=(p_1\hd p_r)$ of $d$,
we write $l(\pi)=r$, $|\pi|=d$, $[\pi]$ is the
irreducible $\FS_d$-module associated to $\pi$, and  $S_{\pi}V$ is the associated irreducible
$GL(V)$ module.
  $Sym(V)$ denotes the symmetric algebra.
  For a variety $X\subset \BP V$,
we let $\hat X\subset V$ denote the corresponding cone
in $V$.
$A_j$ is a vector space of dimension $a_j$ and we
assume $a_j\geq 2$ to
avoid trivialities.
We often write $\s_r=\s_r(\BP A_1^*\ctimes \BP A_n^*)$.

\smallskip

Our main results are as follows:

\begin{theorem}\label{p1pb-1pc-1}
The varieties
$\s_r(\pp 1\times \pp{b-1}\times \pp{c-1})=
\s_r(\BP A^*\times \BP B^*\times \BP C^*)$
are normal, with rational singularities. Their ideal is generated
in degree $r+1$ by the
irreducible modules in the two   flattenings:
$$
\La{r+1}(A\ot B)\ot \La{r+1}C,{\rm \ and \ }
\La{r+1}(A\ot C)\ot \La{r+1}B\subset S^{r+1}(A\ot B\ot C)
$$
The redundancy in the above description is the irreducible module
$S^{r+1}A\ot\La{r+1}B\ot\La{r+1}C$.
\end{theorem}

\begin{theorem}\label{4factor} 
Fix positive integers $a,b,c,d$.
The variety $\s_2(\pp {a-1}\times \pp {b-1}\times \pp {c-1}\times \pp {d-1})=\s_2(\BP A^*\times\BP B^*\times \BP C^*\times \BP D^*)$
is
arithmetically Cohen-Macaulay. Its ideal is generated in degree
three
by the modules defining the subspace variety $Sub_{2222}$
(namely $\La 3 A\ot \La 3(B\ot C\ot D)$ plus permutations minus
redundancies)
and two copies of  the module
$S_{21}A\ot S_{21}B\ot S_{21}C\ot S_{21}D$ 
which   arise from
the flattenings of the form   $(A\ot B)\ot (C\ot D)$.
\end{theorem}

Note that {\it a priori} there are three
modules obtained 
from flattenings
 but
they  only span two independent
  copies of
$S_{21}A\ot S_{21}B\ot S_{21}C\ot S_{21}D$,   see
Equation \eqref{4factorex} and     
Remark \ref{fourfactorexplain}  below. 
The first set of modules
in Theorem \ref{4factor} may be thought of as arising from
the flattenings of the form $A\ot (B\ot C\ot D)$.

The assertion regarding the generators of the ideal
is the Garcia-Stillman-Sturmfels conjecture for four
factors \cite{GSS}, discussed further in \S\ref{genredsect}.

\begin{theorem}\label{3factor} 
Fix positive integers $a,b,c\geq 3$.
The variety $\s_3(\pp {a-1}\times \pp {b-1}\times \pp {c-1})=\s_3(\BP A^*\times\BP B^*\times \BP C^*)$
is
arithmetically Cohen-Macaulay. Its ideal is generated in degree
four
by    the module
$S_{211}A\ot S_{211}B\ot S_{211}C$ 
which   arises from Strassen's commutation condition.
\end{theorem}

{\it A priori} there are three 
modules 
obtained
by flattenings but
they can only span the unique copy of $S_{211}A\ot S_{211}B\ot S_{211}C$ 
 in $S^4(A\ot B\ot C)$, see Equation \eqref{211mult}.

The equations 
 arising from Strassen's commutation condition 
originated in \cite{Strassen}. A discussion of them
in language compatible with this paper can be
found in \cite{LMsecb}.

\medskip

Remarkably, in each of these cases, the
ideal is generated in the minimal possible
degree ($k+1$ for $\s_k$, see \cite{LMsec}).

\subsection*{Overview}

In \S\ref{subpfsect} we prove all the necessary facts about subspace varieties and we deduce Theorem \ref{p1pb-1pc-1}.

In   \S\ref{genredsect} we describe   Garcia-Stillman-Sturmfels conjecture from \cite{GSS}, and a reduction of
it (Theorem \ref{cmenoughforgss}).

The remainder of the proofs   proceed in two steps. First,
in \S\ref{inheritsect}
we show that  the generators of the ideal
of    secant varieties of Segre
varieties 
 can be deduced from solving  the
basic cases of $\s_r(\BP^{r-1}\ctimes \BP^{r-1})$ (Proposition \ref{gensofsr}),
and moreover the
arithmetically Cohen Macaulay  (ACM) property   holds
in any given case if it holds for the
relevant basic
case  plus   a technical hypothesis on modules occurring
in the minimal free resolution of the ideal in the basic cases
(Lemma \ref{acmlemma}).
To prove the ACM property is inherited we 
use a relative version of the machinery of \cite{weyman}.
Namely, inside the desingularization of the subspace variety, we consider
a  subbundle that gives a partial
desingularization of $\s_r$ and whose fibers are isomorphic
to   the basic case,
  and push down
the minimal free resolution of this   subbundle.
Then we study the \lq\lq relative version\rq\rq\  of this resolution on the desingularization of $Sub_{r\hd r}$.
Our results follow from the analysis of the terms of this complex
of sheaves. 
The methods from \cite{weyman} allow 
us to establish two key facts (Lemma \ref{nohighercoh}). First, the higher cohomology of the terms of this complex vanishes. 
Second,   the sections of the terms
 are maximal Cohen-Macaulay modules supported in $Sub_{r\hd r}$.
(The proof of this second fact is the most subtle point
in this paper.)  Lemma \ref{nohighercoh} allows
us to compute the length of a minimal free resolution
of $\s_r$ under certain assumptions
described in  Lemma \ref{acmlemma}.


The basic cases for Theorems \ref{4factor} and \ref{3factor} 
are $\s_2(\pp 1\times \pp 1\times \pp 1\times \pp 1)$ and
$\s_3(\pp 2\times \pp 2\times \pp 2)$.
Respectively in  \S\ref{4p1sect} and  \S\ref{3p2sect}  
we prove these varieties
are ACM, determine generators of their ideals,
and show the technical hypotheses necessary
to apply Lemma \ref {acmlemma} hold.
 Unfortunately this step
utilizes a computer calculation.

It is interesting to ask if the the ACM property holds
for general secant varieties of Segre varieties. 
From our approach it follows that the ACM property for   
$\sigma_r(\BP A_1^*\ctimes \BP A_n^*)$, with $dim A_j\ge r$ would follow from checking the assumptions of Lemma \ref {acmlemma} for the variety $\s_r(\pp {r-1}\ctimes \pp {r-1})$ ($n$-factors).

Since we use results from representation theory, commutative
algebra, and the geometric method of \cite{weyman} throughout,
we begin in \S\ref{review} with brief remarks
from these areas.

\subsection*{Acknowledgment} We thank the anonymous referee
for very useful suggestions to improve the exposition of this paper.

\section{Review from representation theory, commutative
algebra and the geometric method of \cite{weyman}}\label{review}

\subsection{Syzygies}
We summarize from \cite{weyman} (5.1.1-3,5.4.1):
\begin{theorem}\cite{weyman}\label{weymanthm}
 Let $Y\subset \BP V$ be a
variety and suppose there is a
projective variety $B$ and a vector bundle
$E\ra B$ that is a subbundle of a trivial bundle
$\underline V \ra B$ with $\underline V_z\simeq V$ for
$z\in B$ such that $  E\ra \hat Y$ is a desingularization.  Write $\eta=E^*$ and $\xi=(\underline V/E)^* $

If
the sheaf cohomology groups
$H^i(B,S^d\eta)$ are all zero for $i>0$
and if the linear maps
$H^0(B,S^d\eta)\ot V^*\ra H^0(B,S^{d+1}\eta)$
are surjective for all $d\geq 0$, then
\begin{enumerate}
\item
$\hat Y$ is normal, with rational singularities 

\item The coordinate ring $K[\hat Y]$ satisfies
  $K[\hat Y]_d\simeq H^0(B,S^d\eta)$.

\item The 
vector space of minimal generators of the ideal of $\hat Y$ in
degree $d$ is isomorphic to  
$H^{d}(B,\La{d+1}\xi)$.

\item If moreover $Y$ is a $G$-variety
and the desingularization is $G$-equivariant,
then the identifications above are as $G$-modules.
\end{enumerate}
\end{theorem}

More generally, in the situation of Theorem \ref{weymanthm},
$\oplus_j H^j(\La{i+j}\xi)$ is isomorphic to
the $i$-th term in the minimal free resolution
of $Y$, and even a \lq\lq twisted\rq\rq\ version
of this result holds which we recall and explain
when it is used in \S\ref{inheritsect}.

\subsection{Representation theory}
Let $V=A_1\otc A_n$. Let $G=GL(A_1)\ctimes GL(A_n)$.
The varieties $\s_r$ are $G$-varieties so we
should study their ideals as $G$-modules. The first
step in doing this is to
  decompose
$S^dV$ into $G$-isotypic components.
Recall that to a partition $\pi$ we associate a representation 
$[\pi]$ 
of the symmetric group on $d$ letters $\FS_d$ and a representation 
$S_{\pi}W$
of the general linear group $GL(W)$. Both  groups act on
$W^{\ot d}$ and
  each group is the commuting subgroup  of the other.
The $GL(W)$-isotypic decomposition of $W^{\ot d}$ is
$W^{\ot d}=\oplus_{|\pi|=d}[\pi]\ot S_{\pi}W$.

\begin{proposition}(\cite{LMsec}, 4.1) The  $G=GL(A_1)\ctimes GL(A_n)$
isotypic decomposition of $S^d(A_1\otc A_n)$ is 
$$S^d(A_1\otc A_n)= \bigoplus_{|\pi_1|=\cdots =|\pi_k|=d}([\pi_1]\otc [\pi_n])^{\FS_d} \ot
S_{\pi_1}A_1\otc S_{\pi_k}A_k,$$
where $([\pi_1]\otc [\pi_k])^{\FS_d}$ denotes the space of $\FS_d$-invariants
(i.e., instances of the trivial representation of $\FS_d$) in $[\pi_1]\ot \cdots\ot [\pi_n]$. \end{proposition}

\smallskip

Note in particular that the decomposition of $S^d (A_1\ot\cdots\ot A_n )$ is uniform, i.e. if  $\tdim A_i\geq l(\pi_i)$ (so the corresponding module is non-zero), then the multiplicity   $([\pi_1]\otc [\pi_k])^{\FS_d}$ does not depend on  the $\tdim A_i$. 

\smallskip

The multiplicity of $S_{\pi_1}A_1\otc S_{\pi_k}A_k$
in $S^d(A_1\otc A_n)$, which is $\tdim ([\pi_1]\otc [\pi_k])^{\FS_d}$,  
can be computed using   characters in low degrees,
although there is no general closed form formula.
Let
$\chi_{\pi_j}: \FS_d\ra\BC$ denote the   character  of $[\pi_j]$, then 
$$
\tdim ([\pi_1]\otc [\pi_n])^{\FS_d}  =\frac 1{d!}
\sum_{\a\in\FS_d} \chi_{\pi_1}(\s) 
  \cdots\chi_{\pi_n}(\s)
$$ 
(see, e.g., \cite{macdonald}).
 
For example:

\begin{equation}\label{2facdecomp}
([\pi_1]\ot  [\pi_2])^{\FS_d}=\delta_{\pi_1,\pi_2}
\end{equation}  i.e. only symmetric terms occur with multiplicity one,

\begin{equation}\label{4factorex}
\tdim ([(2,1)],[(2,1)],[(2,1)],[(2,1)])^{\FS_3}=2.
\end{equation}

and 

\begin{equation}\label{211mult}
\tdim ([(2,1,1)],[(2,1,1)],[(2,1,1)])^{\FS_4}=1.
\end{equation}

\begin{remark}\label{fourfactorexplain} Assume $A_1 ,A_2 ,A_3 ,A_4$ have all dimension 2. Then   $\La3 (A_i\ot A_j )= S_{2,1}A_i\ot S_{2,1}A_j$. Thus any two flattenings $\La 3 (A_i\ot A_j)\ot\La 3 (A_k\ot A_l)$, embedding this representation into $S^3 (A_1\ot A_2\ot A_3\ot A_4 )$ via $3\times 3$ minors of a $4\times 4$ matrix  span the isotypic component
of
$S_{21}A_1\ot S_{21}A_2\ot S_{21}A_3\ot S_{21}A_4$
in $S^3(A_1\otc A_4)$.
\end{remark}

\subsection{Commutative algebra}

Let $V$ be a $K$-vector space, let $A=Sym(V)$, which we
consider as the algebra of polynomials on $V^*$.
For a graded $A$-module $M$, 
$pd_A(M)$, the {\it projective dimension} of $M$, denotes the length
of a minimal free resolution of $M$ as an $A$-module.

For a homogeneous ideal $I\subset A$, we let
$Z_I\subset V^*$ denote its associated variety (the
zero set of the polynomials in $I$). Similarly, the
{\it support} of an $A/I$-module is
$Z_{Ann(M)}\subset Z_I\subset V^*$.  

\begin{definition} $A/I$ is a 
{\it Cohen-Macaulay ring} iff
$$pd_A (A/I)= \tcodim(Z_I, V^*).$$
An $A/I$-module $M$ is a {\it maximal Cohen-Macaulay module} iff
$$pd_A (M)= \tcodim(Z_I, V^*).$$
An affine variety $Z\subset V^*$ is
{\it arithmetically Cohen-Macaulay (ACM)} if its coordinate ring $K[Z]$
is a Cohen-Macaulay ring, i.e.,  the length of a minimal free resolution
of $K[Z]$ as an $A$-module
equals the codimension of $Z$. 
\end{definition}

The following classical result   follows, e.g., from
\cite{eisenbud}, Theorem 18.15.a.

\begin{theorem}\label{eisquote}Notations as above.
Let   $I\subset A$ be a homogeneous ideal,
let $Z=Z_I\subset V^*$   and let $Z_{sing}$ be
its singular locus.  Assume $A/I$ is Cohen-Macaulay,
then $A/I$ is reduced iff $\tcodim(Z_{sing},Z)\geq 1$.
\end{theorem}

We also note the
following standard Commutative Algebra result,  which essentially
says
that a generically reduced irreducible algebraic variety has an non-empty open subset of smooth points.

\begin{proposition} \label{acmgenred}
 If an affine variety $Z\subset V$ is generically reduced, then $\tcodim(Z_{sing},Z)\geq 1$.
 \end{proposition}

\section{The subspace varieties and their defining ideals}\label{subpfsect}

\begin{theorem}\label{subthm} The subspace varieties
$Sub_{b_1\hd b_n}$
are normal, with rational
singularities. Their ideal is generated in degrees
$b_{j}+1$ for $1\leq j\leq n$ by the irreducible modules
in 
$$\La{b_j+1}A_j\ot\La{b_j+1}(A_1\otc   A_{j-1}\ot A_{j+1} \otc A_n),
$$
such that  (reordering such that $b_1\leq b_2\leq\cdots \leq b_n$)
the partitions $S_{\pi_i}A_i$ that occur for $i\leq j$ have $l(\pi_i)\leq b_i$,
unless $b_i=b_j$, in which case we also allow $l(\pi_i)=b_i+1$.

In particular, if all the $b_i=r$, the ideal of $Sub_{r\hd r}$ is generated
in degree $r+1$ by the irreducible modules appearing
in  
$$\La{r+1}A_j\ot \La{r+1}(A_1\otc A_{j-1}\ot A_{j+1}\otc A_n)
$$
for $1\leq j\leq n$
(minus   redundancies).
\end{theorem}

\begin{proof}
 
First note that the ideal of
$Sub_{b_1\hd b_n}$ consists of all modules
$S_{\pi_1}A_1\otc S_{\pi_n}A_n$ occurring in $S^d(A_1\otc A_n)$
where each $\pi_j$ is a partition of $d$ and at least one $\pi_j$
has $l(\pi_j)>b_j$. Also, notice, that for each $j$ the ideal consisting of representations 
$S_{\pi_1}A_1\otc S_{\pi_n}A_n$ occurring in $S^d(A_1\otc A_n)$
where  $l(\pi_j)>b_j$ is generated in degree $b_j +1$ by
$$\La{b_j+1}A_j\ot\La{b_j+1}(A_1\otc   A_{j-1}\ot A_{j+1} \otc A_n),
$$
because it is just the ideal for rank at most $b_j$ tensors in
the tensor product of two vector spaces.
After reordering of summands  so $b_1\le\ldots\le b_n$ an  elementary induction by degree completes the argument regarding
generators of the ideal.

\smallskip

To prove the results on the singularities,
consider the product of Grassmannians 
$$B=G(b_1,A_1^*)\ctimes G(b_n,A_n^*)
$$ 
and the  
bundle 
\begin{equation}\label{defE}
p:  \cR_{1}\otc \cR_n\ra B
\end{equation}
 where $\cR_j$ is the tautological
rank $b_j$ subspace bundle over $G(b_j,A_j^*)$.
Then the total space $\tilde Z$ of $\cR_{1}\otc \cR_n$ maps to $ A_1^*\otc A_n^*$.
We let $q:\tilde Z\ra  A_1^*\otc A_n^*$ denote this map which
 gives a desingularization of $Sub_{b_1\hd b_n}$.
(A general element of $Sub_{b_1\hd b_n}$ is of the
form $[a^1_1\otc a^1_n +\cdots +a^{b_n}_1\otc a^{b_n}_n]$
where $\tdim\langle a^1_j\hd a^{b_n}_j\rangle = b_j$, so
it has a unique preimage under $q$.)

By Theorem \ref{weymanthm}.1,
with $\eta=(\cR_{1}\otc \cR_n)^*$, we need to show

\smallskip

(i.) $H^i(B,S^d\eta)=0$ for all $i>0$, for all $d\geq 0$

\smallskip
(ii.) $H^0(B,S^d\eta)\ot (A_1\otc A_n)
\ra H^0(B,S^{d+1}\eta)$ is surjective for all $d\geq 0$.

\smallskip

To see (i.) holds, note that $\eta=\cR_{1}^*\otc \cR_n^*$, and thus $S^d(\eta)$, is
 homogeneous, completely reducible,
and the factors are tensor products on Schur functors on $\cR_i^*$.
Each of these irreducible
factors is ample (in fact, a quotient bundle of
a trivial bundle) thus     the Bott-Borel-Weil Theorem implies
$S^d(\eta)$ has no higher cohomology (cohomology of an
irreducible bundle can occur at most in one degree). 

To see (ii.), 
  the ring of sections of $Sym (\eta)$ is generated in degree 0 because
the description of the ideal of $Sub_{b_1\hd b_n}$ given above shows that,   the multiplication
  map is induced  by the multiplication in   $Sym (A_1\ot\cdots\ot A_n )$ after mod-ing out
the span of the representations $S_{\pi_1}A_1\otc
S_{\pi_n}A_n$ satisfying   $l(\pi_j)> b_j$ for some $j$. But the Littlewood-Richardson rule (e.g. \cite{weyman}, Theorem (2.3.4)) implies that in the tensor product of two representations $S_{\pi_1} V\otimes S_{\pi_2} V$ we have only the representations $S_{\pi_3} V$ with 
the Young diagram of $\pi_3$ containing both 
the diagrams of $\pi_1$ and $\pi_2$
as sub-diagrams, so if a representation $S_{\pi_1}A_1\otc S_{\pi_n}A_n$ satisfies $l(\pi_j)\le b_j$ for all $j$, and it appears
in $(S_{\mu_1}A_1\otc S_{\mu_n}A_n)\ot (A_1\otc A_n)$
then $l(\mu_j)\leq b_j$ for all $j$ as well.
\end{proof}

\begin{proof}[proof of Theorem \ref{p1pb-1pc-1}]
 Theorem \ref{subthm} and
Strassen's result \cite{Strassen}   that $\s_r(\pp 1\times\pp {r-1}
\times \pp{r-1})=\BP (K^2\ot K^r\ot K^r)$
(which is easily established using Terracini's lemma)   imply
$\BP (Sub_{2,r,r})=\s_r(\BP A^*\times \BP B^*\times \BP C^*)$ when $b,c\geq r$.
(One always has $\s_r(\BP A^*\times \BP B^*\times \BP C^*)\subseteq
\BP (Sub_{r,r,r})$  and Strassen's result establishes
the reverse inclusion.)
\end{proof}

\section{The varieties $Flat^{\overline a}_r$ and the GSS conjecture}\label{genredsect}

A variant on the subspace varieties is as follows.
Let $\overline a=(a_1\hd a_n)$ and define 
$I_{Flat^{\overline a}_r}$ to be the ideal
generated by the modules $\La {r+1}A_I\ot \La {r+1}A_J\subset
S^{r+1}(A_1\otc A_n)$ as $I,J$ range over complementary
subsets of 
  $\{1\hd n\}$.
We let $Flat^{\overline a}_r$ denote the corresponding variety.
Just as with subspace varieties, we have
$\s_r(\BP A_1^*\ctimes \BP A_n^*)\subseteq Flat^{\overline a}_r$.

Garcia, Stillmann and Sturmfels \cite{GSS} conjectured that
$I_{Flat^{\overline a}_2}=I_{\s_2(\BP A_1^*\ctimes \BP A_n^*)}$.
We refer to this statement as to the {\it GSS conjecture}.
In \cite{LMsec} the conjecture was proven when 
$\overline a=(a_1,a_2,a_3)$, and moreover it was shown that
as sets, $Flat^{\overline a}_2 = \s_2(\BP A_1^*\ctimes \BP A_n^*) $
for all $n$.
Since $\s_2(\BP A_1^*\ctimes \BP A_n^*) $ is reduced and irreducible,
and $Flat^{\overline a}_2 $ is irreducible, to
prove the conjecture  it would be sufficient to show $Flat^{\overline a}_2 $
is reduced.

The application to the GSS conjecture is

\begin{theorem}\label{cmenoughforgss}
If $Flat^{\overline a}_2 $ is arithmetically Cohen-Macaulay, then the GSS conjecture holds.
\end{theorem}

\begin{proof}

We first show

\begin{proposition}\label{genred} 
$Flat^{\overline a}_2 $
is generically reduced.
\end{proposition}

\begin{proof} Fix bases $(a^{s_i}_i)$ in each $A_i$
and let $\phi_{j_1\hd j_n}$ be linear coordinate functions
on $A_1^*\otc A_n^*$. A general element of $Flat^{\overline a}_2$ is of
the form $x=a^1_1\otc a^1_n + a^2_1\otc a^2_n$, i.e., it
has coordinates $\phi_{1\hd 1}=\phi_{2\hd 2}=1$ and all other
coordinates zero (this, and the assertion about the codimension
follows by using the identification 
as sets of $Flat^{\overline a}_2$ with $\hat\s_2$). We show that
at $x$, the differentials
of a set of generators of $I_{Flat^{\overline a}_r}$ span a subspace
of $T^*_x(A_1^*\otc A_n^*)$ equal to the codimension
of $Flat^{\overline a}_2$.
In algebraic language, we show that the localization
of $Sym(A_1\otc A_n)/I_{Flat^{\overline a}_2}$ at $x$
has codimension equal to $\tcodim (Flat^{\overline a}_2)$.
  $T^*_x\s_2$  is spanned by
$d\phi_{j_1\hd j_n}|_x$ where $n-1$ of the $ j_1\hd j_n $ are
neither $1$ nor $2$. Fix some $p<n$ and consider
the $(a_1\cdots a_p)\times (a_{p+1}\cdots a_n)$ matrix
  corresponding to the flattening $(A_1\otc A_p)\ot (A_{p+1}\otc A_n)$.
Examining the differentials of its three by three minors
at $x$, all are zero except the differentials of minors 
containing $\phi_{1\hd 1}$ and $\phi_{2\hd 2}$, which will
have a unique nonzero term $d\phi_{i_1\hd i_n}|_x$.
For any splitting we recover all the 
$d\phi_{i_1\hd i_n}|_x$ where none of the $i_s$ are $1$ or $2$.
In general we recover all the $d\phi_{i_1\hd i_n}|_x$ that are
neither in the row or column containing $\phi_{1\hd 1}$ or
$\phi_{2\hd 2}$. Thus
if we want a term with $k$ indices equal to $1$ and $l$ indices
equal to $2$, then (ignoring order for the moment)
as long as $k<n-p$ and $l<p$ there is clearly no problem.
To get a different order, just permute the factors.

\end{proof}

To conclude the proof of Theorem \ref{cmenoughforgss} we use Theorem \ref{eisquote} and Proposition \ref{acmgenred}. 

\end{proof}

\begin{example} Consider the case $n=4$ and each $a_i=2$. Here are matrices
respectively  for
the splittings $(A_1\ot A_2)\ot (A_3\ot A_4)$
and $(A_1\ot A_3)\ot (A_2\ot A_4)$.

$$
 \left(\begin{matrix} \phi_{1,1,1,1}&\phi_{1,2,1,1}
&\phi_{2,1,1,1}&\phi_{2,2,1,1}
\\
\phi_{1,1,1,2}&\phi_{1,2,1,2}&\phi_{2,1,1,2}&\phi_{2,2,1,2}
\\
  \phi_{1,1,2,1}&\phi_{1,2,2,1}&\phi_{2,1,2,1}&\phi_{2,2,2,1}
\\
\phi_{1,1,2,2}&\phi_{1,2,2,2}&\phi_{2,1,2,2}&\phi_{2,2,2,2}
\end{matrix} \right)$$
 
$$  
\left(\begin{matrix} \phi_{1,1,1,1}&\phi_{1,1,2,1}
&\phi_{2,1,1,1}&\phi_{2,1,2,1}\\
\phi_{1,1,1,2}&\phi_{1,1,2,2}&\phi_{2,1,1,2}&\phi_{2,1,2,2}\\
  \phi_{1,2,1,1}&\phi_{1,2,2,1}&\phi_{2,2,1,1}&\phi_{2,2,2,1}\\
\phi_{1,2,1,2}&\phi_{1,2,2,2}&\phi_{2,2,1,2}&\phi_{2,2,2,2}
\end{matrix} \right)$$

The $d\phi_{ijkl}|_x$ where $\{  i,j,k,l\}=\{1,1,2,2\}$
  each appear in the differentials of the eight relevant
(i.e., those containing both $\phi_{1111}$ and $\phi_{2222}$)  $3\times 3$ minors.
\end{example}

 We resolve 
the four factor case of the GSS conjecture  as a consequence of
Lemma \ref{acmlemma} and  Proposition \ref{4p1case}.

\section{Ideals and the ACM property are inherited}\label{inheritsect}

\begin{definition}
Given
vector spaces $A_j'\subset A_j$
and  a module $S_{\pi_1}A_1'\otc S_{\pi_n}A_n'
\subset S^d(A_1'\otc A_n')$,
we say the module $S_{\pi_1}A_1\otc S_{\pi_n}A_n$
correspondingly realized as a submodule of $S^d(A_1\otc A_n)$
is {\it inherited} from $S_{\pi_1}A_1'\otc S_{\pi_n}A_n'$.

Note that if $S_{\pi_1}A_1'\otc S_{\pi_n}A_n'$ is nonzero,
we have
$S_{\pi_1}A_1'\otc S_{\pi_n}A_n'\subset I(\s_r(Seg(\BP A_1'
\ctimes \BP A_n'))$
iff
$S_{\pi_1}A_1 \otc S_{\pi_n}A_n \subset I(\s_r(Seg(\BP A_1 
\ctimes \BP A_n)) $. This property is called
{\it inheritance} in \cite{LMsec}.

\end{definition}
\begin{proposition}\label{gensofsr}
Let $\tdim A_1\hd \tdim A_n\geq r$.
The generators of the ideal of
$\s_r(\BP A_1^*\ctimes \BP A_n^*)$ are given by the
modules generating the ideal of
$Sub_{r\hd r}$ and the modules inherited from the
modules generating the ideal of
$\s_r(\pp{r-1}\ctimes \pp {r-1})$ ($n$-factors).
\end{proposition}
\begin{proof}
The irreducible modules generating the
ideal of $Sub_{r\hd r}$ are all in degree $r+1$ and
are the irreducible submodules of $\La{r+1}A_j\ot \La{r+1}
(A_1\otc A_{j-1}\ot A_{j+1}\otc A_n)$, so in particular they
all contain a partition with   $r+1$ parts.
The irreducible modules generating the ideal of 
$\s_r(\pp{r-1}\ctimes \pp {r-1})$ cannot contain
a partition with more than $r$ parts.

Now say some module 
$S_{\pi_1}A_1\otc S_{\pi_n}A_n$ is in $I(\s_r(Seg(\BP A_1 
\ctimes \BP A_n))$.
We must show it is generated from our candidate generators.
 If any
$\pi_j$ has more than $r$ parts, then it is already in
the ideal generated by $Sub_{r\hd r}$ so we are done.
But now if each $\pi_j$ has length at most $r$, then
the same module must also be in the ideal
of $\s_r(\pp {r-1}\ctimes \pp {r-1})$. 
\end{proof}

Over the Grassmannian $G(r,A_j^*)$, we let
$\cR_j$, $\cQ_j$ respectively denote the
rank $r$ (resp. rank $a_j-r$) tautological
subspace (resp. quotent) vector bundles.
Recall the bundle
$\eta=\cR_1^*\otc \cR_n^*$. Let
$\cB=Sym(\eta)$.

\begin{lemma}\label{nohighercoh} 
Let $\pi_j=(p_{j,1}\hd p_{j,r})$
 be partitions.  Consider the sheaf
$${\cM}:= \otimes_{j=1}^n S_{\pi_j }{\cR}_j^*\otimes {\cB}.$$
  
\begin{enumerate}
\item Assume that $p_{j,1}\ge -a_j +1$ for $1\le j\le n$. Then     ${\cM} $ is acyclic.
\item Assume that $ p_{j,1}\ge 0$ and $p_{j,1}\le r^{n-1}-r$ for $1\le j\le n$. Then the $Sym(A_1\otc A_n)$-module  $H^0(B,{\cM})$, which is
supported in $Sub_{r\hd r}$,
is a maximal Cohen-Macaulay module.
\end{enumerate}
\end{lemma} 

\begin{proof}
The first assertion is a straightforward application of the
Bott-Borel-Weil theorem.
 
The second assertion is the most subtle point of
this paper. To prove it, we   use the duality theorem
\cite{weyman}, Theorem 5.1.4, which we now recall.

For any vector bundle $\cV\ra B$,
   following \cite{weyman}, Theorem 5.1.4, define
the twisted dual vector bundle
$$
\check {\cV} =K_B\ot \La {\trank \xi}\xi^*
\ot \cV^*
$$
where $\cV^* $ denotes the ordinary dual vector bundle, $K_B$ is the canonical bundle of $B$, and $\xi=(A_1^*\otc A_n^*\ot  \cO_B/
\cR_1\otc \cR_n)^*$.
Then \cite{weyman}, Theorem 5.1.4,
asserts that
\begin{equation}\label{duality}
F(\check\cV)_j=F(\cV)^*_{j+\tdim B -\trank\xi }.
\end{equation}

\smallskip
 
We claim that 
under the hypotheses of the lemma, the rightmost nonzero term in $F(\check\cM)_{\bullet}$ is the 
zero-th.
To see this note that $K_{G(r,a)}=S_{a -r\hd a -r}\cR\ot S_{r\hd r}\cQ^*$,
which up to tensoring with a trivial bundle
(powers of the bundle $(\La{a_1}A_1\otc \La{a_n}A_n)\ot \cO_B$)
 is isomorphic to
$S_{a \hd a }\cR$, and, up to tensoring with a trivial bundle, 
$$
\Lambda^{\trank \xi}\xi^* \simeq S_{r^{n-1}\hd r^{n-1}}\cR_1^*\otc S_{r^{n-1}\hd r^{n-1}}\cR_n^* 
$$ 

Write $\pi_i=(p_{i,1}\hd p_{i,r})$.
So up to tensoring with a trivial line bundle,
$$
\check \cM\simeq
S_{(r^{n-1}-a_1-p_{1,r})\hd (r^{n-1}-a_1-p_{1,1})}\cR_1^*\otc 
S_{(r^{n-1}-a_n-p_{n,r})\hd (r^{n-1}-a_n -p_{n,1})}\cR_n^*
$$
Thus, if for each $i$ we have $p_{i,1}\le r^{n-1}-r$, then,
applying (1), $\check\cM\otimes \cB$ has no higher cohomology, and the complexes $F(\cM)$ and $F(\check\cM)$ have length equal to the codimension of the subspace variety $Sub_{r\hd r}$ which equals $(\trank \xi-\tdim B)$.
\end{proof}

\begin{lemma}\label{acmlemma}
If $\s_r(\pp {r-1}\ctimes \pp {r-1})$ ($n$-factors,
with $n\geq 3$) is arithmetically Cohen-Macaulay
with the property that no module occurring in its minimal
free resolution contains a partition whose first part is
greater than $r^{n-1}-r$,   then
$\s_r(\BP A_1^*\ctimes \BP A_n^*)$ is arithmetically Cohen-Macaulay
when $\tdim A_i\geq r$ for $1\leq i\leq n$,.
\end{lemma}

\begin{proof} Notatations as above.
Consider  the desingularization of the subspace
variety $Sub_{r\hd r}$ and the resulting
vector bundle $E=\cR_{1}\otc   \cR_{n}$ as in 
Equation \eqref{defE} with each $b_j=r$. Each fiber
$(\cR_{1}\otc \cR_{n})_x$ of $E$
over $x\in B=G(r,A_1^*)\ctimes   G(r,A_n^*)$
is just $\BC^r\otc   \BC^r$ and
we may consider the subvariety $Z\subset \tilde Z$
such that $  Z_x= \hat\s_r(\BP  (\cR_{1})_x\ctimes   \BP(\cR_{n})_x)$.
(Recall that
  $\tilde Z$ is
the total space of the bundle $E$.) 
$Z$ gives a partial desingularization of $\hat\s_r$.

Under our hypotheses, there is a minimal free resolution $G_{\bullet}$
of $Z_x$
where $G_0=\bold A:=  Sym(A_1\otc A_n) $, $G_1$ is a sum of modules
$S_{\pi_1}A_1\otc S_{\pi_n}A_n\ot \bold A(-k)$ where $k\geq r+1$, 
and the length of the resolution of $G_{\bullet}$ is
the codimension of $Z_x$ in $\BP E_x$, namely
$L:=r^n-r^2n+r(n-1)$, as $\s_r(\pp{r-1}\ctimes \pp{r-1})$ is
of the expected dimension $(rn+1)(r-1)$ as long as $n\geq 3$.

By \cite{weyman} Proposition (5.1.1), part b), $\cB=Sym(\eta)$ is a sheaf of algebras isomorphic
to $p_*(\cO_{\tilde Z})$.

We form a complex of sheaves of $\cB$-modules from
$G_{\bullet}$ by replacing each $G_i$ with the
sheaf $\cG_i$ obtained by replacing the 
Schur functors of the vector spaces $A_1\hd A_n$ with the corresponding
Schur functors of the sheaves $\cR_{1}^*\hd \cR_{n}^*$.

We have projections $q:Z\ra \hat\s_r$ and $p: Z\ra B$.
We   have $p_*(\cO_Z)= \cB/d(\cG_1 )$  as $d(\cG_1 )$ is
the subsheaf of $\cB$ consisting of the local functions
on $\tilde Z$ that vanish on $Z$.

Our complex of sheaves of $\cB$-modules $\cG_{\bullet}$
is such that each term is a sum of terms of the form
$$
S_{\pi_1}\cR_{1}^*\otc S_{\pi_n}\cR_{n}^*\ot \cB.
$$
Each term is homogeneous and completely reducible, with each irreducible
summand having nonzero $H^0$,
so in particular no term has any higher cohomology.
 
Define a complex $M_{\bullet}$ of $\bold A$-modules
by letting $M_j:= H^0(B,\cG_j)$.

The minimal free resolution of the ideal of $\s_r$
is the minimal resolution of the cokernel of the complex $M_{\bullet}$.
Indeed, by Proposition \ref{gensofsr}, the cokernel $M_0/\tim (M_1)$
is exactly   $K[\s_r]$ because $M_0$ consists
of functions on the subspace variety and $M_1$ the ideal of
the secant variety inside the subspace variety.


\smallskip

To obtain a
not necessarily minimal resolution of the cokernel
$K[\s_r]$  of the map $M_1\ra M_0$, 
one can proceed by iterating the mapping cone construction as follows.
Let $F_{j\bullet}$ be a resolution of $M_j$ for each $j$.
We obtain a double complex, the tail of which is
$$
\begin{array}{cccccc}
\downarrow   &  &\downarrow &  & &  \\
  F_{L,1} & \ra &   F_{L-1,1} & \ra  &  &   \\
 \downarrow   &  &\downarrow &  & &  \\
  F_{L,0} & \ra &   F_{L-1,0} & \ra  &  &   \\
\downarrow & &\downarrow & & &   \\
   M_L & \ra&M_{L-1}  & \ra &  & 
\end{array}
$$
We replace this tail by using the mapping cone
construction (e.g. \cite{eisenbud}), where
we replace $F_{L-1,j}$ by modules $\tilde F_{L-1,j}=F_{L,j-1}\op F_{L-1,j}$,
and $\tilde F_{L-1,\bullet}$ becomes the last column 
of the new  complex.  We iterate this procedure
until we end up with
a picture
$$
\begin{array}{cccccc}
\downarrow   &  &\downarrow &  & &  \\
\tilde F_{11} & \ra &   F_{01} &   &  &   \\
 \downarrow   &  &\downarrow &  & &  \\
\tilde F_{10} & \ra &   F_{00} &   &  &   \\
\downarrow & &\downarrow & & &   \\
 \tilde M_1 & \ra&M_0  & \ra &\cC &\ra 0
\end{array}
$$
where the $F_{0\bullet}$ is a resolution of $M_0$
and $\tilde M_1$ is the term replacing $M_1$ after
having iterated the mapping cone construction, and
$\tilde F_{1\bullet}$ its resolution. 

The final product of this procedure is a
possibly nonminimal resolution $\tilde F_{0\bullet}$ of 
$K[\s_r]$, whose $j$-th term is
$ \tilde F_{0j} =\oplus_{a+b=j} F_{a,b}$.

\smallskip

But by Lemma \ref{nohighercoh}, 
the modules
$M_i$ are maximal Cohen-Macaulay, hence the
lengths of their minimal free resolutions all
equal  $\tcodim Sub_{r\hd r}=(\trank \xi-\tdim B)$.

 
But now the complexes $F(\cG_i)_{\bullet}$ give the resolutions
of the $M_i$, so when we apply the iterated cone construction,
the longest possible length of the possibly nonminimal resolution
of $K[\s_r]$ is 
\begin{align*}
 \trank \xi&-\tdim B +\tcodim\s_r(\pp {r-1}\ctimes \pp{r-1}) \\
&=
 (a_1\cdots a_n-r^n)-r(a_1+\cdots +a_n-nr) +(r^n-r^2n+r(n-1))\\
&=a_1\cdots a_n-r(a_1+\cdots +a_n)+rn-r,
\end{align*}
but 
$$\tcodim \s_r(\BP A_1^*\ctimes \BP A_n^*)
=a_1\cdots a_n-1 - [r(a_1+\cdots + a_n -n)+(r-1)].
$$
We see that the (possibly non-minimal) resolution is of minimal length
and that  length equals the codimension of $\s_r(\BP A_1^*\ctimes \BP A_n^*)$, hence 
$\s_r(\BP A_1^*\ctimes \BP A_n^*)$ is Cohen-Macaulay.
\end{proof}

\section{Case of $\s_2(\pp 1\times \pp 1\times\pp 1\times\pp 1)$}\label{4p1sect}

\begin{proposition}\label{4p1case}
The variety $\s_2(\pp {1}\times \pp {1}\times \pp {1}\times \pp {1})=\s_2(\BP A^*\times\BP B^*\times \BP C^*\times \BP D^*)$
is arithmetically Cohen-Macaulay. Its ideal is generated in degree
three
by two copies of  the module
$S_{21}A\ot S_{21}B\ot S_{21}C\ot S_{21}D$ 
which   arise from
the flattenings of the form   $(A\ot B)\ot (C\ot D)$.
\end{proposition}
\begin{proof}
Let $\bold A=Sym(A\ot B\ot C\ot D)$ and let
$I$ denote the ideal generated by
the relevant two copies $S_{21}A\ot S_{21}B\ot S_{21}C\ot S_{21}D$ (see Remark \ref{fourfactorexplain}). We thank Anurag
Singh for calculating the minimal free resolution
of $\bold A/I$, which we denote $G_\bullet$ with   terms as follows
$$G_\bullet : 0\rightarrow \bold A^{12}(-10)\rightarrow \bold A^{48}(-9)\rightarrow \bold A^{57}(-8)\rightarrow \bold A^{20}(-6)\oplus \bold A^{48}(-5)\rightarrow$$
$$\rightarrow \bold A^{78}(-4)\rightarrow \bold A^{32}(-3)\rightarrow \bold A.$$
Note that $6=\tcodim \s_2(\pp 1\times\pp 1\times\pp 1\times \pp 1)=
\tcodim Flat^{2222}_2$ and
since this coincides with the length of the minimal free resolution
we conclude  that $\bold A/I$ is Cohen-Macaulay. But  we know that $Flat^{2222}_2$ is  Cohen-Macaulay, so by Theorem \ref{cmenoughforgss} the GSS conjecture follows in this case.
\end{proof}

Theorem \ref{4factor} follows because if
we express the resolution in terms of modules,
each module
$$S_{\pi_1}A\ot S_{\pi_2}B\ot S_{\pi_3}C\ot S_{\pi_4}D
$$
that
occurs in some $G_j$  indeed satisfies the property that the first part
of each $\pi_i$ is less or equal to $6$.
This can be calculated directly by examining the maps produced
by Macaulay2 and then finding the equivariant form of the
resolution explicitly (which we reproduce below).
To see it more directly,
note that since the coordinate ring is Cohen-Macaulay, the dual of this resolution
is also an acyclic complex. This means that every representation $S_{\pi_1}A\ot S_{\pi_2}B\ot S_{\pi_3}C\ot S_{\pi_4}D$
appearing in the resolution has to have partitions $\pi_1$, $\pi_2$, $\pi_3$, $\pi_4$ that are contained in partitions of some representation occurring in the top of the resolution. But in the top piece, dimension considerations show immediately that the partitions are  
$\FS[ (6,4)(5,5)(5,5)(5,5) ]$
and thus all the partitions $\pi_j$ have all parts $\le 6$ as required.

\begin{remark} The resolution $G_{\bullet}$
expressed as a direct sum of 
  $GL(A)\times GL(B)\times GL(C)\times GL(D)$-modules
  is as follows.
Denote the $i$-th term in this resolution by $G_i$.  Let $(a,b)(c,d)(e,f)(g,h)$ denote    $S_{(a,b)}A\otimes S_{(c,d)}B\otimes S_{(e,f)}C\otimes S_{(g,h)}D$. The terms in the resolution $G_\bullet$ have to be symmetric under permuting the spaces   so 
we let  $\FS [(a,b)(c,d)(e,f)(g,h)]$  denote the direct sum of all
distinct
 tensor products of Schur functors which are obtained from a given one by permutations of $A,B,C,D$.
\begin{align*}
 G_6 &=\FS[ (6,4)(5,5)(5,5)(5,5) ]\otimes \bold A(-10),\\
 G_5 &= \FS [(5,4)(5,4)(5,4)(5,4)]\otimes \bold A^3 (-9  ),\\
 G_4 &=\FS [(5,3)(5,3)(4,4)(4,4)]\otimes \bold A(-8)\oplus \FS [(4,4)(4,4)(4,4)(4,4)]\otimes \bold A^3(-8),\\
 G_3 &= \FS [(3,2)(3,2)(3,2)(3,2)]\otimes \bold A^3 (-5)\oplus  \FS [(5,1)(3,3)(3,3)(3,3)]\otimes \bold A(-6),\\
 G_2 &= \FS [(3,1)(2,2)(2,2)(2,2)]\otimes \bold A^2 (-4)\oplus \FS [(3,1)(3,1)(2,2)(2,2)]\otimes \bold A(-4),\\
 G_1 &= \FS [(2,1)(2,1)(2,1)(2,1)]\otimes \bold A^2 (-3),\\
 G_0 &=\bold A.
\end{align*}
\end{remark}

\section{Case of $\s_3(\pp 2\times \pp 2\times \pp 2)$}\label{3p2sect}

\begin{proposition}\label{3p2case}
The variety $\s_3(\pp {2}\times \pp {2}\times \pp {2} )=\s_3(\BP A^*\times\BP B^*\times \BP C^* )$
is arithmetically Cohen-Macaulay. Its ideal is generated in degree
four
by   the module
$S_{211}A\ot S_{211}B\ot S_{211}C $ 
which   arises from Strassen's commutation condition.
(A priori there are three  copies obtained this way but
they are all isomorphic submodules of $S^4(A\ot B\ot C)$.)
\end{proposition}

\begin{remark} In \cite{GSS} it is stated without proof that Strassen's equations
  generate the ideal of $\s_3(\pp {2}\times \pp {2}\times \pp {2} )$.
In personal communication,
the authors informed us that their assertion comes from a Macaulay2 calculation similar to
the one Hal Schenck calculated for us below.
\end{remark}

\begin{proof} Let $\bold A= Sym(A\ot B\ot C)$ and let $I$
denote the ideal generated by Strassen's polynomials.
We thank Hal Schenck for calculating the minimal
free resolution of $\bold A/I$ (using Macaulay2 \cite{mac2}) which we denote
$G_{\bullet}$ with terms as follows:
$$G_\bullet : 0\rightarrow \bold A(-15)\oplus \bold A^{30}(-12)\rightarrow \bold A^{189}(-11)\rightarrow \bold A^{351}(-10)\rightarrow \bold A^{223}(-9)\oplus \bold A(-6)\rightarrow$$
$$\rightarrow \bold A^{30}(-6)\oplus \bold A^{27}(-5)\rightarrow \bold A^{27}(-4)\rightarrow \bold A.$$

Since $\tcodim \s_3(\pp 2\times \pp 2\times \pp 2)= 6$, the ideal is ACM. If we show that it is generically reduced, then by applying Theorem \ref{eisquote} we see it is reduced. But we know by \cite{Strassen} that it defines the secant variety set-theoretically, so it has to be the defining ideal.

Let $a_1,a_2,a_3$ be a basis of $A^*$ and similarly for
$B,C$. Let $\phi_{ijk}$ denote coordinates on $A^*\ot B^*\ot C^*$
with respect to the bases and let $x=[a_1\ot b_1\ot c_1+a_2\ot b_2\ot c_2+a_3\ot b_3\ot c_3]\in \s_3$ which corresponds to
the point with $\phi_{iii}=1$ and all other coordinates zero.
The conormal space to
$\s_3$ at $x$ is given by the span of $d\phi_{ijk}|_x$ with
$i,j,k$ distinct. We must show that the differentials of
the polynomials coming from Strassen's equations give all
of these. By symmetry it is sufficient to show they
give one of these.

Following \cite{LMsecb},   write
$$T=a_1\ot X + a_2\ot Y\ot a_3\ot Z$$
 where
  $X,Y,Z$
are represented as
$3\times 3$ matrices with respect to   bases of $B^*,C^*$,
then 
$$
P_{1st}(T)=
\sum_{j,k}(-1)^{j+k}(\tdet X^{\hat \j}_{\hat k})
(Y^j_tZ^s_k-Y^s_kZ^j_t)
$$
where $X^{\hat \j}_{\hat k}$ is $X$ with its $j$-th row and
$k$-th column removed. The polynomials
$P_{2st},P_{3st}$ are obtained by exchanging the roles
of $X$ respectively with $Y,Z$.

For example, $dP_{1st}|_x=0$ unless $s=2,t=3$ or $s=3,t=2$. In
the first case one obtains $dP_{123}|x=-d\phi_{123}$ and
in the second one obtains $dP_{132}|_x=d\phi_{132}$.
In general, for distinct $i,j,k$,  $dP_{ijk}|_x=\pm d\phi_{ijk}$
and it is zero otherwise.
\end{proof}

Theorem \ref{3factor} follows by observing that 
expressing the resolution in terms of modules,
if
$S_{\pi_1}A\ot S_{\pi_2}B\ot S_{\pi_3}C$
occurs in some $G_j$, that the first part
of each $\pi_i$ is less than $6$. As mentioned
above, this can be read off of the equivariant
form of the resolution (which can be deduced
from the Macaulay2 printout) or determined by dimension
considerations or more directly using the top piece
of the resolution
$$G_6=\FS[(5,5,5)(5,5,5)(5,5,5)]\oplus \FS[(6,3,3)(4,4,4)(4,4,4)]$$
and noting that here
all the partitions   have all parts $\le 6$ as required,
as in the proof of Theorem \ref{4factor}.

\begin{remark}\label{lastrem}
The equivariant form of the resolution $G_\bullet$
is as follows: Let $(a,b,c)(d,e,f)(g,h,i)$ 
denote the tensor product  $S_{(a,b,c)}A\otimes S_{(d,e,f)}B\otimes S_{(g,h,i)}C$. The terms in the resolution $G_\bullet$ have to be symmetric under permuting the spaces $A,B,C$ so we denote by $\FS [(a,b,c)(d,e,f)(g,h,i)]$ the direct sum of all 
distinct tensor products of Schur functors which are obtained from a given one by permutations of $A,B,C$.  The terms are as follows:
\begin{align*}
G_6&=\FS[(5,5,5)(5,5,5)(5,5,5)]\oplus \FS[(6,3,3)(4,4,4)(4,4,4)],\\ 
 G_5& = \FS[(5,3,3)(4,4,3)(4,4,3)]\oplus \FS[(4,4,3)(4,4,3)(4,4,3)],\\ 
G_4& = \FS[(4,3,3)(4,4,2)(4,4,2)]\oplus \FS[(4,3,3)(4,3,3)(4,3,3)], \\
 G_3& = \FS[(4,4,1)(3,3,3)(3,3,3)]\oplus \FS[(4,3,2)(4,3,2)(3,3,3)]\oplus \FS[(3,3,3)(3,3,3)(3,3,3)]
\\
 & \ \ \ \oplus\FS[(2,2,2)(2,2,2)(2,2,2)],\\
 G_2 &= \FS[(4,1,1)(2,2,2)(2,2,2)]\oplus \FS[(2,2,1)(2,2,1)(2,2,1)],\\
 G_1& = \FS [(2,1,1)(2,1,1)(2,1,1)],\\
 G_0& =\FS[(0,0,0)(0,0,0)(0,0,0)].
\end{align*}
This can be deduced from the Macaulay printout by analyzing the weights of basis elements in the resolution.
\end{remark}

\vspace{1cm}

{\small
\noindent {\sc Joseph M. Landsberg}, 
Department of Mathematics,
  Texas A\&M University,
  Mailstop 3368,
  College Station, TX 77843-3368, USA

\smallskip
 
\noindent {\sc Jerzy Weyman}, 
Department of Mathematics,
  Northeastern University,
  360 Huntington Avenue,
Boston, MA 02115, USA 
}

\smallskip

\end{document}